\documentclass[11pt,a4paper]{article}

\usepackage[T1]{fontenc}
\usepackage[utf8]{inputenc}
\usepackage{lmodern}
\usepackage[margin=1.12in]{geometry}
\usepackage{microtype}
\usepackage{amsmath,amssymb,amsthm,mathtools,mathrsfs}
\usepackage{tikz-cd}
\usepackage{enumitem}
\usepackage{aliascnt}
\usepackage{needspace}
\usepackage[numbers,sort&compress]{natbib}
\usepackage[unicode,hidelinks,hypertexnames=false]{hyperref}
\usepackage[capitalise,noabbrev]{cleveref}

\tikzcdset{
  every label/.append style={font=\small,inner sep=2.2pt},
  cells={nodes={inner sep=3pt}},
  arrows={line width=0.45pt},
  column sep=4.4em,
  row sep=3.2em,
  proof square/.style={column sep=5.2em,row sep=3.4em},
  proof strip/.style={column sep=4.0em,row sep=3.3em}
}
\setlist{itemsep=0.25em,topsep=0.5em}
\numberwithin{equation}{section}
\allowdisplaybreaks[1]
\hypersetup{
  pdftitle={Absolute descent and the Benabou--Roubaud theorem},
  pdfauthor={Fernando Lucatelli Nunes},
  pdfsubject={Left decalage, absolute descent, and monadicity},
  pdfkeywords={absolute descent, left decalage pseudomonad, successor contraction, pseudo-Kan extension}
}

\newtheorem{theorem}{Theorem}[section]
\newaliascnt{proposition}{theorem}
\newtheorem{proposition}[proposition]{Proposition}
\aliascntresetthe{proposition}
\crefname{proposition}{Proposition}{Propositions}
\Crefname{proposition}{Proposition}{Propositions}
\newaliascnt{lemma}{theorem}
\newtheorem{lemma}[lemma]{Lemma}
\aliascntresetthe{lemma}
\crefname{lemma}{Lemma}{Lemmas}
\Crefname{lemma}{Lemma}{Lemmas}
\newaliascnt{corollary}{theorem}
\newtheorem{corollary}[corollary]{Corollary}
\aliascntresetthe{corollary}
\crefname{corollary}{Corollary}{Corollaries}
\Crefname{corollary}{Corollary}{Corollaries}
\theoremstyle{definition}
\newaliascnt{definition}{theorem}
\newtheorem{definition}[definition]{Definition}
\aliascntresetthe{definition}
\crefname{definition}{Definition}{Definitions}
\Crefname{definition}{Definition}{Definitions}
\newaliascnt{example}{theorem}
\newtheorem{example}[example]{Example}
\aliascntresetthe{example}
\crefname{example}{Example}{Examples}
\Crefname{example}{Example}{Examples}
\theoremstyle{remark}
\newaliascnt{remark}{theorem}

\aliascntresetthe{remark}
\crefname{remark}{Remark}{Remarks}
\Crefname{remark}{Remark}{Remarks}

\AddToHook{env/theorem/before}{\Needspace{6\baselineskip}}
\AddToHook{env/proposition/before}{\Needspace{6\baselineskip}}
\AddToHook{env/lemma/before}{\Needspace{5\baselineskip}}
\AddToHook{env/corollary/before}{\Needspace{5\baselineskip}}
\AddToHook{env/definition/before}{\Needspace{5\baselineskip}}
\AddToHook{cmd/subsection/before}{\Needspace{7\baselineskip}}
\newcommand{\op}{\mathrm{op}}
\newcommand{\id}{\mathrm{id}}
\newcommand{\Hcat}{\mathcal H}
\newcommand{\Kcat}{\mathcal K}
\newcommand{\C}{\mathcal C}
\newcommand{\Lx}{\mathcal L}
\newcommand{\Cat}{\mathbf{Cat}}

\newcommand{\two}{\mathbf 2}
\newcommand{\PsRan}{\operatorname{PsRan}}
\newcommand{\Desc}{\operatorname{Desc}}

\newcommand{\suc}{\operatorname{suc}}
\newcommand{\Dmon}{\mathbb D}
\newcommand{\Smon}{\mathbb S}

\newcommand{\Contr}{\operatorname{Contr}}
\newcommand{\PsAlg}{\operatorname{PsAlg}}
\newcommand{\Adj}{\mathsf{Adj}}

\title{Absolute descent and the
B\'enabou--Roubaud theorem}
\author{Fernando Lucatelli Nunes\\[0.4em]
  \small Centre for Mathematics of the University of Coimbra (CMUC),\\[-0.1em]
  \small Department of Mathematics, University of Coimbra, 3000--143 Coimbra, Portugal\\[-0.1em]
  \small\href{mailto:fln@uc.pt}{\texttt{fln@uc.pt}}}
\date{}

\begin{document}

\maketitle
\begin{abstract}
	We continue the programme developed in \emph{Pseudo-Kan extensions and
		descent theory} by studying descent for augmented pseudocosimplicial
	objects in their own right. Let \(\mathcal H\) be a \(2\)-category,
	and let \(\suc=\mathsf{1}+(-)\) be the successor monad on the augmented
	simplex category. Precomposition with \(\suc\) induces the
	\textit{d\'ecalage} pseudomonad \(\Smon\) on the \(2\)-category of augmented
	pseudocosimplicial objects with values in \(\mathcal H\). We prove
	that every successor contraction makes the corresponding augmented
	pseudocosimplicial object effective for descent. This structure
	is preserved by postcomposition with every pseudofunctor.
	Consequently, the resulting descent object is an absolute bilimit.
	
	A first application concerns split epimorphisms. We recover the
	Janelidze--Tholen characterization of absolute effective descent
	morphisms given in \emph{Facets of Descent II}. Moreover, we prove
	that one direction holds for every pseudofunctor
	\(A\colon C^{\op}\to\Hcat\), where \(\Hcat\) is an arbitrary
	\(2\)-category.
	
	A further consequence of our main theorem is that every
	\(\Smon\)-pseudoalgebra is of absolute descent. In particular,
	this holds for every free \(\Smon\)-pseudoalgebra, and hence for
	every one-sided d\'ecalage. This makes the shift-effectiveness
	hypothesis in the B\'enabou--Roubaud theorem of
	\emph{Pseudo-Kan extensions and descent theory} automatic. Thus, under
	the adjointness and Beck--Chevalley hypotheses, the descent and
	Eilenberg--Moore factorizations are equivalent whenever they exist.
	We thus obtain a stronger version of the B\'enabou--Roubaud theorem for
	arbitrary augmented pseudocosimplicial objects, which we call the
	formal B\'enabou--Roubaud theorem. Applied, in particular, to indexed
	\v{C}ech objects associated with pseudofunctors
	\(A\colon C^{\op}\to\Hcat\), this theorem extends the classical
	B\'enabou--Roubaud theorem to an arbitrary \(2\)-category \(\Hcat\).
\end{abstract}


\Needspace{12\baselineskip}
\section{Introduction}
\label{sec:introduction}
Descent theory studies the image of a functor, typically one induced
by a change of base. In this interpretation, objects in the source
play the role of global objects, while their images in the target play
the role of local realizations. The descent problem is to recognize,
from within the target, which local objects arise from global ones.

A bare description of the essential image rarely provides an
intrinsic criterion for recognizing its objects. Instead, one equips
objects in the target with coherent compatibility data and proves that
they come from the source precisely when they admit such data. One
must also identify the morphisms that preserve this structure, thereby
reconstructing the source, up to equivalence, as a category of local
objects equipped with descent data. Descent thus transforms an
external problem about the image of a functor into an internal
recognition principle.

This pattern occurs throughout mathematics. For modules, bundles,
sheaves, and indexed categories, a global object determines local data
together with canonical compatibility information. Effective descent
provides the converse: every coherent system of such data arises,
essentially uniquely, from a global object
\cite{Grothendieck1960,Grothendieck1971,Giraud1964,JanelidzeTholen1994}.

The programme developed in \cite{Nunes2018} places Grothendieck
descent theory in a general two-dimensional setting, organized around
two complementary principles:
\begin{enumerate}
  \item \emph{Bilimit viewpoint.} In an ambient \(2\)-category,
  descent objects are conical bilimits of pseudocosimplicial objects;
  see \cite[Theorem~4.11]{Nunes2018}.

  \item \emph{Pseudomonadic viewpoint.} For pseudomonadic
  pseudofunctors, essential images can be recognized through the
  existence of suitable pseudoalgebra structures. Idempotent
  pseudomonads are particularly well suited to this purpose: their
  pseudoalgebra structures are essentially unique and may be
  characterized by pseudoretractions of the unit; see
  \cite[Section~2]{Nunes2018}.
\end{enumerate}

\paragraph{Bilimit viewpoint.}
The first principle can be expressed in more general terms: descent is
expressed by the appropriate notion of limit of a cosimplicial object.
It therefore makes sense to develop descent for such objects in their
own right, independently of the change-of-base constructions from
which they may arise.

In dimension \(1\), the conical limit of a cosimplicial object in an
ordinary category is the equalizer of its first two cofaces, whenever
either exists. This is natural:
compatibility in dimension \(1\) is expressed by equality, and the
cosimplicial identities extend any morphism equalizing these two
cofaces uniquely to a cone over the entire diagram. For example,
the sheaf condition identifies sections over a covered open set
with the equalizer of the two maps restricting families of local
sections to pairwise overlaps: compatible local sections glue uniquely.

In dimension \(2\), in the pseudo
setting, compatibility is expressed by isomorphisms satisfying
coherence conditions. Here, the equalizing equation is replaced by a
specified invertible \(2\)-cell satisfying normalization and a cocycle
condition. These are precisely the conditions defining Grothendieck
descent \cite{Grothendieck1971,JanelidzeTholen1994,Pavlovic1991}.
Indeed, by \cite[Theorem~4.11 and Remark~4.4]{Nunes2018}, the usual
descent object, defined by these data and conditions in the first
three nonempty degrees, is equivalent to the conical bilimit of the
entire pseudocosimplicial object.

More explicitly, let \(\dot\Delta\) be the category of finite ordinals
and monotone maps, including the empty ordinal \(\mathsf0\), and let
\[
j\colon\Delta\hookrightarrow\dot\Delta
\]
be the full inclusion of the nonempty finite ordinals. We regard each ordinal
\[
\mathsf n=\{0\leq\ldots\leq n-1\}
\]
as a category. An augmented pseudocosimplicial object in \(\Hcat\)
is a pseudofunctor \(A\colon\dot\Delta\longrightarrow\Hcat\);
see \cite[Definition~2.1]{Nunes2016} for the definition of a
pseudofunctor.
We denote by
\begin{equation}
d_n^i\colon A(\mathsf n)\longrightarrow A(\mathsf{n+1}),
\qquad s_n^i\colon A(\mathsf{n+1})\longrightarrow A(\mathsf n)
\end{equation}
the coface induced by the monotone injection omitting \(i\), for
\(0\leq i\leq n\), and the codegeneracy induced by the monotone
surjection identifying \(i\) and \(i+1\), for \(0\leq i<n\).

Writing \(d=d_0^0\), the relevant part of \(A\) is
\[
\begin{tikzcd}[column sep=3.8em]
  A(\mathsf 0) \ar[r,"d"]
  & A(\mathsf 1)
    \ar[r,bend left=40,"d_1^0"]
    \ar[r,bend right=40,"d_1^1"']
  & A(\mathsf 2)
    \ar[l,"s_1^0" description]
    \ar[r,bend left=40,"d_2^0"]
    \ar[r,"d_2^1" description]
    \ar[r,bend right=40,"d_2^2"']
  & A(\mathsf 3)  \, \cdots 
\end{tikzcd}
\]
\noindent
where the positive part \(A \circ j\) contains the descent data, while \(A(\mathsf0)\) is seen as 
the proposed reconstruction.

The augmented pseudocosimplicial structure equips \(d\) with a
canonical descent datum for \(A \circ j\). Hence, whenever \(\Desc(A \circ j)\)
exists in \(\Hcat\), its universal property induces the descent
factorization
\begin{equation}
\begin{tikzcd}[column sep=6.0em,row sep=3.2em]
	A(\mathsf 0) \ar[r,"k"] \ar[dr,"d"',""{name=D}]
	& \Desc(A \circ j) \ar[d,"u"]
	\\
	& A(\mathsf 1)
	\ar[phantom,from=1-2,to=D,
	"\cong" description,pos=0.65]
\end{tikzcd}
\end{equation}
given by the canonical invertible \(2\)-cell \(u\circ k\cong d\).
We write \(k=k_A\) and \(u=u_A\) when the diagram needs to be specified.
Here \(u\) is the forgetful morphism, while
\[
k\colon A(\mathsf 0)\longrightarrow\Desc(A\circ j)
\]
is the descent comparison. The augmented pseudocosimplicial object
\(A\) is of \emph{effective descent} precisely when \(k\) is an
equivalence.

\paragraph{Pseudomonadic viewpoint.}
The second principle returns descent to an image problem, now in the
\(2\)-category
\[
[\dot\Delta,\Hcat]_{PS}
\]
of augmented pseudocosimplicial objects, pseudonatural transformations,
and modifications. The proposed reconstruction \(A(\mathsf 0)\) can be replaced
by the universal one, \(\Desc(A\circ j)\), producing an augmented object of
effective descent. In \cite{Nunes2018}, this completion is formalized
by an idempotent pseudomonad constructed from pointwise pseudo-Kan
extensions.

More precisely, suppose that the relevant conical bilimits exist for
all diagrams, so that restriction along \(j\) has a pointwise right
pseudo-Kan-extension biadjoint
\[
j^*\colon[\dot\Delta,\Hcat]_{PS}
\rightleftarrows[\Delta,\Hcat]_{PS}:\PsRan_j.
\]
The induced \emph{descent pseudomonad} is
\[
\Dmon=\langle j\rangle=\PsRan_j\circ j^*.
\]
Up to equivalence, \(\Dmon\) leaves the positive part of \(A\)
unchanged and replaces \(A(\mathsf 0)\) by \(\Desc(A\circ j)\). Since \(j\) is fully
faithful, \(\Dmon\) is idempotent
\cite[Theorem~3.30]{Nunes2018}. Its unit
\[
\eta_A\colon A\longrightarrow\Dmon(A)
\]
is the descent comparison: its degree-zero component is \(k_A\), while
its positive-degree components are equivalences, and may be taken to be isomorphisms. Consequently,
\[
\begin{aligned}
	A\text{ is of effective descent}
	&\quad\Longleftrightarrow\quad
	\eta_A\colon A\longrightarrow\Dmon(A)
	\text{ is an equivalence}
	\\
	&\quad\Longleftrightarrow\quad
	A\text{ admits a }\Dmon\text{-pseudoalgebra structure}.
\end{aligned}
\]
The effective-descent objects are therefore precisely the objects in
the essential image of the forgetful pseudofunctor
\[
\PsAlg(\Dmon)\longrightarrow[\dot\Delta,\Hcat]_{PS},
\]
or, equivalently, the objects fixed up to equivalence by descent
completion. The pseudomonadic viewpoint thus provides an organizing
principle for the bilimit viewpoint: it expresses the same
reconstruction in terms of pseudoalgebra structure, while bringing
effective descent back to the problem of characterizing the essential
image of a forgetful pseudofunctor.

This viewpoint does not require descent objects for every
\(\Hcat\)-valued diagram. We may restrict to the full
\(2\)-subcategory of augmented pseudocosimplicial objects \(A\)
for which \(\Desc(A\circ j)\) exists. Descent completion preserves this
subcategory and defines an idempotent pseudomonad on it, with the
same characterization of effective descent by pseudoalgebra
structure. We make this restriction precise in
\cref{prop:restricted-descent-completion}. For Theorem~A below,
even existence of the particular descent object will follow from
the successor contraction.

\subsection*{Absolute descent}

The present paper asks a stronger stability question. Effective
descent asks whether reconstruction succeeds; absolute descent asks
whether the same reconstruction survives every change of ambient
\(2\)-category. More precisely, an augmented pseudocosimplicial
object \(A\colon\dot\Delta\to\Hcat\) is of \textit{absolute effective descent}
if \(P\circ A\) is of effective descent for every pseudofunctor
\(P\colon\Hcat\to\Kcat\).

The ordinary case suggests what such a criterion should look like.
A splitting expresses the universal property of an equalizer by
equations preserved by every functor, and hence makes the equalizer
absolute. Since ordinary descent reduces to an equalizer, the problem
is to identify the coherent two-dimensional structure that plays the
role of this splitting.

The required structure comes from the successor monad
\[
\suc=\mathsf{1}+(-)\colon\dot\Delta\longrightarrow\dot\Delta,
\]
which adjoins a new least element. Viewing ordinals as categories,
this is the restriction to \(\dot\Delta\) of the monad that freely
adjoins an initial object. Its unit
\(\delta\colon\id_{\dot\Delta}\to\suc\) has components
\[
\delta_n\colon\mathsf n\longrightarrow\mathsf1+\mathsf n
\]
given by inclusion into the second summand. Its multiplication has
components
\[
\mu_n\colon\mathsf1+\mathsf1+\mathsf n
\longrightarrow\mathsf1+\mathsf n
\]
that identify the two adjoined least elements and act as the identity
on the original \(\mathsf n\)-summand.

For an augmented pseudocosimplicial object \(A\), the components of
\[
\id_A\ast\delta\colon A\longrightarrow A\circ\suc
\]
are the first cofaces \(d_n^0\). A \emph{successor contraction} on
\(A\) is a pseudonatural retraction of \(\id_A\ast\delta\), namely a
pseudonatural transformation
\[
r\colon A\circ\suc\longrightarrow A
\]
equipped with an invertible modification \(r\circ(\id_A\ast\delta)\cong\id_A\).
The main observation of this paper is the following.

\medskip
\noindent\textbf{Theorem A (Absolute descent).}
\emph{Let \(\Hcat\) be any \(2\)-category. If an augmented
	pseudocosimplicial object \(A\colon\dot\Delta\to\Hcat\) admits a
	successor contraction, then \(A\) is of absolute effective descent.}
\medskip

Precomposition with the successor monad induces the
\textit{d\'ecalage pseudomonad}
\[
\Smon=(-)\circ\suc
\]
on the \(2\)-category \([\dot\Delta,\Hcat]_{PS}\) of augmented
pseudocosimplicial objects. We call \(A\circ\suc\) the
\emph{left d\'ecalage} of \(A\), with \emph{left} referring to the
adjoined least element in the ordinal sum \(\mathsf1+(-)\).
The multiplication \(\id_A\ast\mu\) equips it with the free
\(\Smon\)-pseudoalgebra structure on \(A\) and, in particular,
with a successor contraction.

Theorem~A therefore gives the following consequence.

\medskip
\noindent\textbf{Theorem B (Absolute shifts).}
\emph{For every augmented pseudocosimplicial object \(A\),
	the left d\'ecalage
	\[
	A\circ\suc\colon\dot\Delta\longrightarrow\Hcat
	\]
	is of absolute effective descent.}
\medskip

We now place these results within a general mechanism for contractions
of pointed endopseudofunctors. A successor contraction is precisely a
contraction for the point of \(\Smon\).
The key is a canonical pointed comparison \(\Dmon\to\Smon\),
which transfers this structure to the idempotent descent pseudomonad.
Idempotence then gives effective descent, and preservation of
successor contractions under postcomposition gives absoluteness.

\subsection*{Contractions and pseudoalgebras}

The proof of Theorem~A becomes transparent once the point of a
pseudomonad is separated from its multiplication. A pointed
endopseudofunctor \(\mathbb T=(T,\iota)\) consists of an endopseudofunctor
\(T\) on a $2$-category $\Hcat$ and a pseudonatural transformation
\(\iota\colon\id_{\Hcat} \to T\). A \(\mathbb T\)-contraction on
\(X\) is a pseudoretraction of the point:
\[
T(X)\xrightarrow{x}X,
\qquad
x\circ\iota_X\cong\id_X.
\]
There is no multiplication to respect and therefore no associativity
axiom. We write \(\Contr(\mathbb T)\) for the \(2\)-category of
contractions of \(\mathbb T\), with pseudomorphisms and
compatible \(2\)-cells. For a pseudomonad, this notation refers to its
underlying pointed endopseudofunctor.

Whenever the descent pseudomonad is defined, we construct a
canonical pointed comparison \(\chi\colon\Dmon\to\Smon\), where
\(\iota_A=\id_A\ast\delta\) is the d\'ecalage point.
Restriction along \(\chi\) induces a homomorphism
\[
\chi^*\colon\Contr(\Smon)\longrightarrow\Contr(\Dmon).
\]
It sends a successor contraction \(r\colon\Smon(A)\to A\) to
the descent contraction
\[
\Dmon(A)\xrightarrow{\chi_A}\Smon(A)\xrightarrow{r}A.
\]
The unit cell is the composite
\[
r\circ\chi_A\circ\eta_A
\cong r\circ\iota_A
\cong\id_A.
\]
Since \(\Dmon\) is idempotent, this contraction extends
essentially uniquely to a \(\Dmon\)-pseudoalgebra structure.
The extension also respects pseudomorphisms and \(2\)-cells, so
it gives a biequivalence
\(\Contr(\Dmon)\simeq\PsAlg(\Dmon)\).

It gives the hierarchy
\[
\begin{tikzcd}[column sep=7.0em,row sep=3.2em]
  \PsAlg(\Smon) \ar[r,"U_{\Smon}"]
  & \Contr(\Smon) \ar[d,"\chi^*"]
  \\
  \PsAlg(\Dmon)
  & \Contr(\Dmon) \ar[l,"\text{idempotence}"',"\simeq"]
\end{tikzcd}
\]
The first arrow forgets multiplication coherence. Restriction along
\(\chi\) then turns a successor contraction into a splitting of the
descent unit, and idempotence promotes that splitting essentially
uniquely to a descent pseudoalgebra. Theorem~A begins with
\(\Contr(\Smon)\), rather than with \(\PsAlg(\Smon)\), which is why it
applies to contractions that are not associative.

The transferred pseudoalgebra makes the descent unit an equivalence,
so the original augmented object is of effective descent. Successor
contractions are preserved by postcomposition with every
pseudofunctor; applying the same transfer to the image contraction
therefore gives absoluteness. The existence of the descent objects
needed for this argument, without an ambient completeness
assumption, is established in \cref{thm:absolute-descent}.

\subsection*{The two completions in the B\'enabou--Roubaud theorem}

Absolute descent of successor shifts has a second consequence. It
supplies, automatically and in the stronger absolute form, the
shift-effectiveness hypothesis used in the effective-descent part of
the formal B\'enabou--Roubaud theorem of
\cite[Theorem~8.5]{Nunes2018}. It also makes visible the two
pseudoalgebra structures governing that theorem.

The first coface
\[
d  = d_0^0\colon A(\mathsf 0)\longrightarrow A(\mathsf 1)
\]
admits two natural reconstruction procedures. The
pseudocosimplicial structure gives its descent factorization through
\(\Desc(A\circ j)\). If \(d_0^0\) has a left adjoint, it also induces a monad
\[
T_A=d_0^0\circ(d_0^0)_!
\]
on \(A(\mathsf 1)\), and hence an Eilenberg--Moore factorization through
\(A(\mathsf 1)^{T_A}\). The B\'enabou--Roubaud theorem asserts that, under the
Beck--Chevalley condition, these two reconstructions agree.

Suppose, more generally, that all components of \(\id_A\ast\delta\) have
left adjoints and that the corresponding Beck--Chevalley mates
are invertible. By doctrinal adjunction, these componentwise
adjunctions assemble into an adjunction in
\[
[\dot\Delta,\Hcat]_{PS},
\]
or equivalently into a two-variable pseudofunctor
\[
C_A\colon\Adj\times\dot\Delta\longrightarrow\Hcat,
\]
where \(\Adj\) is the walking adjunction. This diagram carries two
independent reconstruction procedures. Descent completion acts in the
pseudocosimplicial variable, while Eilenberg--Moore completion acts in
the adjunction variable. Both are idempotent pseudo-Kan completions,
and Fubini for pointwise pseudo-Kan extensions identifies the two
orders in which they may be performed, whenever the required objects
exist.

The boundary of the resulting completion square is already effective.
The shifted row \(A\circ\suc\) is of absolute descent by Theorem~B.
For every \(n\geq1\), the right adjoint
\[
d_n^0\colon A(\mathsf n)\longrightarrow A(\mathsf{n+1})
\]
is a pseudosection, with retraction \(s_n^0\), and is therefore
monadic. Thus the shifted row carries the required descent
pseudoalgebra, while the positive-degree adjunctions carry the required
Eilenberg--Moore pseudoalgebras. Fubini fills the only remaining corner
in either order and produces an equivalence
\[
\Desc(A\circ j)\simeq A(\mathsf 1)^{T_A}
\]
compatible with the comparison morphisms from \(A(\mathsf 0)\) and the
forgetful morphisms to \(A(\mathsf 1)\).

Three constructions occur here, but only two are interchanged. The
successor pseudomonad supplies the absolute shifted row. The actual
Fubini comparison is between the idempotent descent and
Eilenberg--Moore completions.

\medskip
\noindent\textbf{Theorem C (Formal B\'enabou--Roubaud).}
\emph{Let \(\Hcat\) be any \(2\)-category and let
	\(A\colon\dot\Delta\to\Hcat\) be augmented pseudocosimplicial. Suppose
	that the components of \(\id_A\ast\delta\) have left adjoints and that
	the corresponding Beck--Chevalley mates are invertible. Whenever the
	descent object and the relevant Eilenberg--Moore object exist, the
	descent factorization of \(d\) is equivalent to its
	Eilenberg--Moore factorization. Hence \(A\) is of effective descent if
	and only if \(d\) is monadic.}
\medskip

The existence of these two particular objects is the only ambient
completeness requirement. Both always exist in \(\Cat\), but the
theorem itself is valid in any \(2\)-category in which they exist.

\subsection*{Indexed \v{C}ech descent}

The preceding results concern pseudocosimplicial objects in an
arbitrary \(2\)-category. Classical descent theory is recovered by
applying them to indexed \v{C}ech objects. Consider a pseudofunctor
\[
\Lx\colon\C^{\mathrm{op}}\longrightarrow\Hcat,
\]
where \(\C\) has pullbacks, and a morphism \(p\colon E\to B\) in \(\C\). If the change-of-base
\(1\)-cells have left adjoints and \(\Lx\) satisfies
Beck--Chevalley for pullback squares, then Theorem~C applied to the
indexed \v{C}ech object \(\Lx\circ D_p\) identifies the descent
factorization of \(\Lx(p)\) with its Eilenberg--Moore factorization,
whenever the two objects exist. Consequently,
\[
p\text{ is of effective }\Lx\text{-descent}
\quad\Longleftrightarrow\quad
\Lx(p)\text{ is monadic}.
\]
This is the two-dimensional B\'enabou--Roubaud theorem for arbitrary
\(\Hcat\); taking \(\Hcat=\Cat\) yields its classical form
\cite{BenabouRoubaud1970,JanelidzeTholen1994}.

The final application uses the absolute-descent theorem directly and
requires neither adjoints nor Beck--Chevalley. If \(p\) has a section, inserting it in the first coordinate contracts the entire
augmented \v{C}ech nerve, rather than merely its kernel-pair fork.
Theorem~A then makes the resulting cone an absolute descent cone for
every indexed pseudofunctor with values in any \(2\)-category.

In the \(\Cat\)-valued setting, every split epimorphism is therefore a
morphism of absolute effective descent. Conversely, the single indexed category
\[
\Lx_B(X)=\two^{\C(B,X)}
\]
detects whether \(p\) has a section. This proves that the morphisms of
absolute effective descent are precisely the split epimorphisms and
recovers the Janelidze--Tholen characterization
\cite[Theorem~3.5]{JanelidzeTholen1997}.

\subsection*{Organization}

\Cref{sec:descent} develops descent both as a bilimit and as an
idempotent pseudo-Kan completion. \Cref{sec:contractions} proves the
transfer principle for contractions.
\Cref{sec:successor} constructs the pointed comparison
\(\Dmon\to\Smon\) and proves absolute descent for successor
contractions and one-sided d\'ecalage. \Cref{sec:BR} presents the
formal B\'enabou--Roubaud theorem as an application of the comparison
results of \cite{Nunes2018} and absolute descent of successor shifts. \Cref{sec:split} applies the results to
indexed \v{C}ech objects and proves the Janelidze--Tholen
characterization. \Cref{sec:codual} gives the co-dual comparison
with co-Eilenberg--Moore objects.

\section{Descent as a bilimit and as an image}
\label{sec:descent}

This section makes precise the two viewpoints announced in the
introduction. We describe descent objects as bilimits and then package effective
descent as pseudoalgebra structure for an idempotent pseudo-Kan
completion.

\subsection{Conventions}

We work relative to a fixed hierarchy of universes and leave changes
of universe implicit; see \cite{Levy2018}. The symbol $\Cat$ denotes
the $2$-category of categories, functors, and natural transformations.
Ordinary categories are regarded as locally discrete $2$-categories.
For $2$-categories $\mathcal J$ and $\Hcat$, we
write $[\mathcal J,\Hcat]_{PS}$ for the $2$-category of
pseudofunctors, pseudonatural transformations, and modifications, and
$[\mathcal J,\Hcat]_{PS}(F,G)$ for its hom-category. In particular, this notation
applies when $\mathcal J$ is an ordinary category.

The compositor and unitor of a pseudofunctor $F$ are oriented as
\[
  F(v)\circ F(u)\Longrightarrow F(v\circ u),
  \qquad
  \id_{F(i)}\Longrightarrow F(\id_i).
\]
For a pseudonatural transformation $\tau\colon F\to G$, the
pseudonaturality constraint associated with $u\colon i\to k$ is
oriented as
\[
  G(u)\circ\tau_i\Longrightarrow\tau_k\circ F(u).
\]
We write $\circ$ for composition of $1$-cells and for vertical
composition of $2$-cells, read from right to left, and $\ast$ for
horizontal composition. In particular, whiskering is written
$\id_g\ast\alpha$ or $\alpha\ast\id_f$. Identities are denoted
by $\id$, with a subscript when needed. Pseudofunctor evaluation is
written with parentheses, as in $F(f)$. Associators, unitors, and
structural isomorphisms are
suppressed when coherence determines them. See
\cite{KellyStreet1974,Lack2010} for background.

A $1$-cell $f\colon X\to Y$ is called a \emph{pseudosection} if it
admits a retraction up to an invertible $2$-cell: there are
$r\colon Y\to X$ and $r\circ f\cong\id_X$. This terminology will be used both
in a $2$-category and in a pseudofunctor $2$-category.

\subsection{Pseudocosimplicial objects and descent data}

Let $\dot\Delta$ be the category of finite ordinals and monotone maps,
including the empty ordinal $\mathsf0$, and let $\Delta$ be its full
subcategory on the nonempty ordinals. We write
\begin{equation}
  j\colon\Delta\hookrightarrow\dot\Delta,
  \qquad
  \mathsf n=\{0,\ldots,n-1\},
  \label{eq:j-inclusion}
\end{equation}
where $j$ is the inclusion. A pseudofunctor
$G\colon\Delta\to\Hcat$ is a \emph{pseudocosimplicial object}; a
pseudofunctor $A\colon\dot\Delta\to\Hcat$ is an
\emph{augmented pseudocosimplicial object}. The coface
$\delta_n^i\colon\mathsf n\to\mathsf{n+1}$ omits $i$, and the
codegeneracy $\sigma_n^i\colon\mathsf{n+1}\to\mathsf n$ identifies
$i$ and $i+1$. Their images are denoted by
\[
  d_n^i=A(\delta_n^i),
  \qquad
  s_n^i=A(\sigma_n^i).
\]
The positive part begins as
\begin{equation}
\begin{tikzcd}[column sep=3.8em]
  A(\mathsf 1)
    \ar[r,bend left=40,"d_1^0"]
    \ar[r,bend right=40,"d_1^1"']
  & A(\mathsf 2)
    \ar[l,"s_1^0" description]
    \ar[r,bend left=40,"d_2^0"]
    \ar[r,"d_2^1" description]
    \ar[r,bend right=40,"d_2^2"']
  & A(\mathsf 3)\,\cdots .
\end{tikzcd}
\label{eq:positive-cosimplicial}
\end{equation}
Only the codegeneracy needed for the normalization condition is
displayed.

For $T\in\Hcat$, a descent datum for $G\colon\Delta\to\Hcat$ with vertex $T$ consists of
a $1$-cell $x\colon T\to G(\mathsf 1)$ and an invertible $2$-cell
\begin{equation}
  \theta\colon d_1^1\circ x\Longrightarrow d_1^0\circ x
  \label{eq:descent-cell}
\end{equation}
satisfying normalization and the cocycle equation. In strict notation
these read
\begin{equation}
  \id_{s_1^0}\ast\theta=\id_x,
  \qquad
  (\id_{d_2^0}\ast\theta)\circ(\id_{d_2^2}\ast\theta)
  =\id_{d_2^1}\ast\theta;
  \label{eq:descent-axioms}
\end{equation}
for a pseudofunctor, the canonical compositor and unitor isomorphisms
are inserted in \eqref{eq:descent-axioms}. A morphism
$(x,\theta)\to(y,\vartheta)$ is a $2$-cell
$\beta\colon x\Rightarrow y$ for which the following square commutes
in $\Hcat(T,G(\mathsf 2))$:
\begin{equation}
\begin{tikzcd}[proof square]
  d_1^1\circ x \ar[r,Rightarrow,"\theta"]
          \ar[d,Rightarrow,"\id_{d_1^1}\ast\beta"']
  & d_1^0\circ x \ar[d,Rightarrow,"\id_{d_1^0}\ast\beta"] \\
  d_1^1\circ y \ar[r,Rightarrow,"\vartheta"'] & d_1^0\circ y.
\end{tikzcd}
\label{eq:descent-morphism}
\end{equation}
These data form a category, denoted
$\mathscr D_G(T)$. Precomposition in $\Hcat$ makes
\begin{equation}
  \mathscr D_G\colon\Hcat^{\mathrm{op}}\longrightarrow\Cat
  \label{eq:descent-data-pseudofunctor}
\end{equation}
a pseudofunctor.

\begin{definition}[Descent object]
\label{def:descent-object}
A \emph{descent object} for
$G\colon\Delta\to\Hcat$ is an object $\Desc(G)$ representing
\eqref{eq:descent-data-pseudofunctor}; thus there are equivalences
\[
  \Hcat(T,\Desc(G))\simeq\mathscr D_G(T)
\]
pseudonatural in $T$.
\end{definition}

Equivalently, $\Desc(G)$ is the conical bilimit of $G$. Its description by descent
data uses only the first three nonempty degrees: degree
$\mathsf1$ contains the underlying object, degree $\mathsf2$ the
descent isomorphism, the first codegeneracy its normalization, and
degree $\mathsf3$ the
cocycle equation \cite[Theorems~4.10--4.11]{Nunes2018}.

For an augmented object $A$, write $d=d_0^0$. Its canonical descent
cell is
\begin{equation}
  \kappa
  :=a_{\delta_1^0,\delta_0^0}^{-1}
     \circ a_{\delta_1^1,\delta_0^0}
  \colon d_1^1\circ d\Longrightarrow d_1^0\circ d,
  \label{eq:canonical-descent-cell}
\end{equation}
where $a_{v,u}\colon A(v)\circ A(u)\Rightarrow A(v\circ u)$ is the
compositor of $A$. The pseudofunctor axioms give the normalization
and cocycle identities for $\kappa$. Whenever $\Desc(A\circ j)$
exists, the canonical comparison induces a functor
\[
  \Hcat(T,k_A)\colon
  \Hcat(T,A(\mathsf0))\longrightarrow\Hcat(T,\Desc(A\circ j)).
\]
Under the equivalence defining $\Desc(A\circ j)$, this sends
\[
  z\colon T\to A(\mathsf 0)
  \quad\longmapsto\quad
  (d\circ z,\kappa\ast\id_z).
\]

\subsection{Pseudo-Kan completion}

The bilimit description is equivalently a pseudo-Kan-extension
description \cite[Corollary~3.27 and Theorem~4.11]{Nunes2018}.
Let $G\colon\Delta\to\Hcat$. When the conical bilimit of $G$ exists,
its right pseudo-Kan extension $\PsRan_j(G)$ exists and is
characterized by equivalences
\[
  [\dot\Delta,\Hcat]_{PS}(X,\PsRan_j(G))
  \simeq[\Delta,\Hcat]_{PS}(X\circ j,G)
\]
pseudonatural in $X\colon\dot\Delta\to\Hcat$. Its values are
\begin{equation}
  (\PsRan_j(G))(\mathsf n)
  \simeq
  \{\dot\Delta(\mathsf n,j(-)),G\}_{\mathrm{bi}}.
  \label{eq:pointwise-ran}
\end{equation}
For $n\geq1$, the weight in \eqref{eq:pointwise-ran} is representable,
so the corresponding value is equivalent to $G(\mathsf n)$. Since
$\mathsf0$ is initial in $\dot\Delta$, the weight at $\mathsf0$ is
constant at the terminal category, and hence
\begin{equation}
  (\PsRan_j(G))(\mathsf0)\simeq\Desc(G).
  \label{eq:descent-as-ran}
\end{equation}
For $A\colon\dot\Delta\to\Hcat$, the unit of the
pseudo-Kan-extension biadjunction has degree-zero component
\begin{equation}
  k_A\colon A(\mathsf 0)\longrightarrow\Desc(A\circ j),
  \label{eq:descent-comparison}
\end{equation}
induced by the canonical cone. Its positive components are
equivalences.

The domain of this completion can be specified without assuming
all descent bilimits in \(\Hcat\).
Let \([\Delta,\Hcat]_{PS}^{\mathrm{desc}}\) be the full
\(2\)-subcategory of \([\Delta,\Hcat]_{PS}\) on the diagrams
\(G\) for which \(\Desc(G)\) exists, and let
\([\dot\Delta,\Hcat]_{PS}^{\mathrm{desc}}\) be the full
\(2\)-subcategory of \([\dot\Delta,\Hcat]_{PS}\) on those \(A\)
with \(A\circ j\in[\Delta,\Hcat]_{PS}^{\mathrm{desc}}\).

\begin{proposition}[The domain of descent completion]
\label{prop:restricted-descent-completion}
Restriction and right pseudo-Kan extension give a biadjunction
\[
  j^*\colon[\dot\Delta,\Hcat]_{PS}^{\mathrm{desc}}
  \rightleftarrows
  [\Delta,\Hcat]_{PS}^{\mathrm{desc}}:\PsRan_j.
\]
Its counit is a pseudonatural equivalence. Consequently,
\begin{equation}
  \Dmon_{\Hcat}=\langle j\rangle_{\Hcat}
  :=\PsRan_j\circ j^*
  \label{eq:descent-pseudomonad}
\end{equation}
is an idempotent pseudomonad on
\([\dot\Delta,\Hcat]_{PS}^{\mathrm{desc}}\)
\cite[Theorem~3.30]{Nunes2018}.
\end{proposition}

We write \(\eta\colon\id\to\Dmon_{\Hcat}\) for the unit and
omit the subscript \(\Hcat\) when the ambient \(2\)-category is
fixed. The full subcategory
\([\dot\Delta,\Hcat]_{PS}^{\mathrm{desc}}\) is the domain of
\(\Dmon_{\Hcat}\) throughout.

\begin{definition}[Effective descent]
\label{def:effective-descent}
An augmented pseudocosimplicial object $A$ is \emph{of effective
$j$-descent} if $\Desc(A\circ j)$ exists and the comparison
$k_A$ in \eqref{eq:descent-comparison} is an equivalence.
\end{definition}

\begin{proposition}[Effective descent as a pseudoalgebra]
\label{prop:effective-descent-pseudoalgebra}
For $A\in[\dot\Delta,\Hcat]_{PS}^{\mathrm{desc}}$, the following
conditions are equivalent:
\begin{enumerate}[label=\textup{(\roman*)},leftmargin=2.7em]
  \item $A$ is of effective $j$-descent;
  \item $\eta_A\colon A\to\Dmon_{\Hcat}(A)$ is an equivalence;
  \item $A$ admits a $\Dmon_{\Hcat}$-pseudoalgebra structure.
\end{enumerate}
Such a pseudoalgebra structure is essentially unique.
\end{proposition}

\begin{proof}
The unit is already an equivalence in every positive degree, so the
first two conditions are equivalent by
\eqref{eq:descent-comparison}. The remaining equivalence is the
standard characterization of pseudoalgebras for an idempotent
pseudomonad \cite[Theorem~2.20]{Nunes2018}.
\end{proof}

Thus the forgetful pseudofunctor
\[
  \PsAlg(\Dmon_{\Hcat})\longrightarrow
  [\dot\Delta,\Hcat]_{PS}^{\mathrm{desc}}
\]
identifies its essential image with the effective-descent objects
in \([\dot\Delta,\Hcat]_{PS}^{\mathrm{desc}}\). This formulation
assumes only the descent objects occurring in that domain. To
prove existence from a successor contraction, we shall apply the
same pseudoalgebra argument representably in \(\Cat\), where
descent completion is defined on every diagram.

\subsection{Absolute descent}

\begin{definition}[Absolute descent]
\label{def:absolute-descent}
An augmented pseudocosimplicial object
$A\colon\dot\Delta\to\Hcat$ is \emph{of absolute effective descent},
or simply \emph{of absolute descent}, if $A$ is of effective
$j$-descent and $F\circ A$ is of effective $j$-descent for every
$2$-category $\Kcat$ and every pseudofunctor $F\colon\Hcat\to\Kcat$.
\end{definition}

\section{Contractions of pointed endopseudofunctors}
\label{sec:contractions}

The descent completion of the previous section is an idempotent
pseudomonad. The left d\'ecalage construction of the next section is also a
pseudomonad, but the absolute-descent theorem does not use its
multiplication. To express exactly what the proof needs, we separate
an endopseudofunctor with a point from a pseudomonad with a
multiplication.

\subsection{Pointed endopseudofunctors}

Let $\Kcat$ be a $2$-category.

\begin{definition}[Pointed endopseudofunctor]
\label{def:pointed-endopseudofunctor}
A \emph{pointed endopseudofunctor} on $\Kcat$ is a pair
\[
  \mathbb T=(T,\iota),
  \qquad
  T\colon\Kcat\to\Kcat,
  \qquad
  \iota\colon\id_{\Kcat}\to T,
\]
where $T$ is a pseudofunctor and $\iota$ is pseudonatural.
\end{definition}

\begin{definition}[Contraction]
\label{def:contraction}
A \emph{$\mathbb T$-contraction} on $X\in\Kcat$ consists of a $1$-cell
\[
  x\colon T(X)\longrightarrow X
\]
and an invertible $2$-cell
\begin{equation}
  \bar x\colon x\circ\iota_X\Longrightarrow\id_X.
  \label{eq:contraction-unit}
\end{equation}
Thus a contraction is a pseudoretraction of the point, with its
invertible retraction cell specified.
\end{definition}

A pseudomorphism of contractions
$(f,\bar f)\colon(X,x,\bar x)\to(Y,y,\bar y)$ consists of a $1$-cell
$f\colon X\to Y$ and an invertible $2$-cell
\[
  \bar f\colon f\circ x\Longrightarrow y\circ T(f)
\]
compatible with the unit cells. Explicitly, if
$\iota_f\colon T(f)\circ\iota_X\Rightarrow\iota_Y\circ f$ is the
pseudonaturality cell of the point, the following diagram commutes:
\[
\begin{tikzcd}[column sep=5.0em,row sep=3.0em]
  f\circ x\circ\iota_X
    \ar[r,Rightarrow,"\bar f\ast\id_{\iota_X}"]
    \ar[d,Rightarrow,"\id_f\ast\bar x"']
  & y\circ T(f)\circ\iota_X
    \ar[d,Rightarrow,"\id_y\ast\iota_f"]
  \\
  f
  & y\circ\iota_Y\circ f
    \ar[l,Rightarrow,"\bar y\ast\id_f"]
\end{tikzcd}
\]
A $2$-cell of contractions
$\beta\colon(f,\bar f)\Rightarrow(g,\bar g)$ is a $2$-cell
$\beta\colon f\Rightarrow g$ for which the following diagram commutes:
\[
\begin{tikzcd}[column sep=5.0em,row sep=3.0em]
  f\circ x
    \ar[r,Rightarrow,"\bar f"]
    \ar[d,Rightarrow,"\beta\ast\id_x"']
  & y\circ T(f)
    \ar[d,Rightarrow,"\id_y\ast T(\beta)"]
  \\
  g\circ x \ar[r,Rightarrow,"\bar g"']
  & y\circ T(g)
\end{tikzcd}
\]
We denote the resulting $2$-category by $\Contr(\mathbb T)$.
For a pseudomonad $\mathbb T$, this notation refers to contractions
of its underlying pointed endopseudofunctor.

If $\mathbb T=(T,\iota,\mu)$ is a pseudomonad, a
$\mathbb T$-pseudoalgebra also has an invertible multiplication cell
\[
  x\circ T(x)\Longrightarrow x\circ\mu_X
\]
satisfying the unit and associativity axioms
\cite{Street1972,BlackwellKellyPower1989}. Forgetting this cell gives
a homomorphism
\begin{equation}
  U_{\mathbb T}\colon
  \PsAlg(\mathbb T)\longrightarrow\Contr(\mathbb T).
  \label{eq:forget-multiplication}
\end{equation}
Throughout, \emph{pseudoalgebra} includes both unit and
multiplication coherence; a contraction requires only the retraction
cell \eqref{eq:contraction-unit}.

\subsection{Restriction along a pointed comparison}

\begin{definition}[Pointed comparison]
\label{def:pointed-comparison}
Let $\mathbb T=(T,\iota)$ and $\mathbb S=(S,\sigma)$ be pointed
endopseudofunctors on $\Kcat$. A \emph{pointed comparison}
\[
  (\chi,\omega)\colon\mathbb T\longrightarrow\mathbb S
\]
consists of a pseudonatural transformation $\chi\colon T\to S$ and
an invertible modification
\begin{equation}
  \omega\colon\chi\circ\iota\Longrightarrow\sigma.
  \label{eq:pointed-comparison-cell}
\end{equation}
\end{definition}

At $X$, Condition \eqref{eq:pointed-comparison-cell} is the coherent
triangle
\begin{equation}
\begin{tikzcd}[column sep=5.5em,row sep=3.0em]
  X \ar[r,"\iota_X"] \ar[dr,"\sigma_X"',""{name=SIGMA}]
  & T(X) \ar[d,"\chi_X"]
  \\
  & S(X)
  \ar[Rightarrow,from=1-2,to=SIGMA,
      shorten <=1.0em,shorten >=0.5em,
      "\omega_X"',"\cong"]
\end{tikzcd}
\label{eq:pointed-triangle}
\end{equation}

\begin{proposition}[Restriction of contractions]
\label{prop:contraction-restriction}
Every pointed comparison
\mbox{$(\chi,\omega)\colon\mathbb T\to\mathbb S$} induces a homomorphism over $\Kcat$
\[
  \chi^*\colon
  \Contr(\mathbb S)\longrightarrow\Contr(\mathbb T).
\]
On objects,
\begin{equation}
  (X,s,\bar s)\longmapsto
  \bigl(X,s\circ\chi_X,
  \bar s\circ(\id_s\ast\omega_X)\bigr).
  \label{eq:contraction-restriction}
\end{equation}
\end{proposition}

\subsection{Pseudoalgebras from contractions}

Let $\Dmon=(D,\eta,\mu)$ be an idempotent pseudomonad on $\Kcat$, so
that $\mu\colon D\circ D\to D$ is a pseudonatural equivalence.
We recall the characterization of its pseudoalgebras
\cite[Theorems~2.18 and~2.20]{Nunes2018}, expressed below in terms
of contractions.

\begin{proposition}[Pseudoalgebras for an idempotent pseudomonad]
\label{prop:idempotent-upgrade}
For $X\in\Kcat$, the following conditions are equivalent:
\begin{enumerate}[label=\textup{(\roman*)},leftmargin=2.7em]
  \item $X$ admits a $\Dmon$-pseudoalgebra structure;
  \item $\eta_X\colon X\to D(X)$ is an equivalence;
  \item $\eta_X$ is a pseudosection;
  \item $X$ admits a $\Dmon$-contraction.
\end{enumerate}
Forgetting multiplication gives a biequivalence over $\Kcat$ whose
inverse adds the uniquely determined multiplication cell and
preserves the contraction data, pseudomorphisms and $2$-cells:
\begin{equation}
  U_{\Dmon}\colon\PsAlg(\Dmon)\longrightarrow\Contr(\Dmon),
  \qquad
  E_{\Dmon}\colon\Contr(\Dmon)\longrightarrow\PsAlg(\Dmon).
  \label{eq:idempotent-forgetful-equivalence}
\end{equation}
\end{proposition}

\begin{theorem}[From contractions to pseudoalgebras]
\label{thm:contraction-transfer}
Let $\Dmon$ be an idempotent pseudomonad and let
$\mathbb S=(S,\iota)$ be a pointed endopseudofunctor on $\Kcat$.
A pointed comparison
\[
  (\chi,\omega)\colon\Dmon\longrightarrow\mathbb S
\]
induces a homomorphism over $\Kcat$
\[
  E_{\Dmon}\circ\chi^*\colon
  \Contr(\mathbb S)\longrightarrow\PsAlg(\Dmon).
\]
For a contraction $(X,s,\bar s)$, the resulting pseudoalgebra action is
\begin{equation}
  D(X)\xrightarrow{\chi_X}S(X)\xrightarrow{s}X,
  \label{eq:transferred-action}
\end{equation}
with unit cell
\[
  s\circ\chi_X\circ\eta_X
  \xRightarrow{\ \id_s\ast\omega_X\ }s\circ\iota_X
  \xRightarrow{\ \bar s\ }\id_X,
\]
and its multiplication cell is uniquely determined by this
contraction. On pseudomorphisms and $2$-cells, the homomorphism
acts by restriction along $\chi$, leaving the underlying $1$-cells
and $2$-cells unchanged.
\end{theorem}

\begin{proof}
Compose the homomorphism $\chi^*$ of
\cref{prop:contraction-restriction} with $E_{\Dmon}$ from
\cref{prop:idempotent-upgrade}. The latter supplies multiplication
coherence and proves its compatibility with every restricted
pseudomorphism and $2$-cell.
\end{proof}

When $\mathbb S$ is a pseudomonad, this gives a diagram of
homomorphisms between the corresponding $2$-categories:
\begin{equation}
\begin{tikzcd}[column sep=7.0em,row sep=3.2em]
  \PsAlg(\mathbb S) \ar[r,"U_{\mathbb S}"]
  & \Contr(\mathbb S) \ar[d,"\chi^*"]
  \\
  \PsAlg(\Dmon)
  & \Contr(\Dmon) \ar[l,"E_{\Dmon}"',"\simeq"]
\end{tikzcd}
\label{eq:contraction-hierarchy}
\end{equation}
Every arrow preserves the carrier. The bottom arrow supplies the
uniquely determined multiplication coherence, on both objects and
morphisms, and preserves the given contraction data.

\section{The left d\'ecalage point and absolute descent}
\label{sec:successor}

We now apply the transfer homomorphism of \cref{sec:contractions} to
descent.

\subsection{The left d\'ecalage pseudomonad}

We define
\begin{equation}
  \suc=\mathsf1+(-)\colon\dot\Delta\longrightarrow\dot\Delta
  \label{eq:successor}
\end{equation}
by adjoining a new least element. Thus
$\suc(\mathsf n)=\mathsf{n+1}$, and for
$u\colon\mathsf m\to\mathsf n$,
\[
  \suc(u)(0)=0,
  \qquad
  \suc(u)(i+1)=u(i)+1.
\]
There is a natural transformation
\begin{equation}
  \delta\colon\id_{\dot\Delta}\longrightarrow\suc,
  \qquad
  \delta_n(i)=i+1.
  \label{eq:successor-unit}
\end{equation}
Its component $\delta_n=\delta_n^0$ is the coface omitting the newly
adjoined least element.

The successor is a monad. Its multiplication
$\mu\colon\suc\circ\suc\to\suc$ has components
\begin{equation}
  \mu_n=\sigma_{n+1}^0\colon
  \mathsf{n+2}\longrightarrow\mathsf{n+1},
  \label{eq:successor-multiplication}
\end{equation}
which identify the two newly adjoined least elements. In coordinates,
$\mu_n(0)=\mu_n(1)=0$ and $\mu_n(i+2)=i+1$.
Precomposition therefore defines the \emph{left d\'ecalage pseudomonad}
\begin{equation}
  \Smon_{\Hcat}:=(-)\circ\suc
  \label{eq:successor-pseudomonad}
\end{equation}
on $[\dot\Delta,\Hcat]_{PS}$. Its unit at $A$ is
\begin{equation}
  \iota_A=\id_A\ast\delta\colon A\longrightarrow A\circ\suc,
  \label{eq:successor-point}
\end{equation}
whose degree-$n$ component is $d_n^0$. The adjective \emph{left}
records that $\suc$ adjoins the new element on the left. We retain the name
\emph{successor monad} for $(\suc,\delta,\mu)$ on
$\dot\Delta$. Thus the successor monad acts on the indexing category,
whereas the d\'ecalage pseudomonad acts on augmented
pseudocosimplicial objects. In particular,
\[
  \Smon_{\Hcat}(A)(\mathsf n)=A(\mathsf{n+1}),
  \qquad
  d_{\Smon(A),n}^{\,i}=d_{A,n+1}^{\,i+1},
  \qquad
  s_{\Smon(A),n}^{\,i}=s_{A,n+1}^{\,i+1}.
\]
The structural comparison cells in these identifications are those
of $A$.

\begin{definition}[Successor contraction]
\label{def:successor-contraction}
A \emph{successor contraction} on
$A\colon\dot\Delta\to\Hcat$ is an $\Smon_{\Hcat}$-contraction.
Explicitly, it consists of a pseudonatural
transformation
\[
  r\colon A\circ\suc\longrightarrow A
\]
and an invertible modification
\begin{equation}
  \varepsilon\colon r\circ(\id_A\ast\delta)\Longrightarrow\id_A.
  \label{eq:successor-contraction}
\end{equation}
\end{definition}

Writing $r_n\colon A(\mathsf{n+1})\to A(\mathsf n)$, the component of
\eqref{eq:successor-contraction} is
\[
  \varepsilon_n\colon r_n\circ d_n^0\Longrightarrow\id_{A(\mathsf n)}.
\]
For every $u\colon\mathsf m\to\mathsf n$, pseudonaturality supplies
an invertible $2$-cell
\begin{equation}
  A(u)\circ r_m\Longrightarrow r_n\circ A(\suc(u)).
  \label{eq:successor-pseudonaturality}
\end{equation}
These cells, together with the modification axiom for
$\varepsilon$, are the higher coherence absent from a splitting of
$d_0^0$ alone.

\subsection{From d\'ecalage structure to descent structure}

We work on the full $2$-subcategory
$[\dot\Delta,\Hcat]_{PS}^{\mathrm{desc}}$ on which the descent
pseudomonad $\Dmon_{\Hcat}$ is defined. The next proposition also
establishes that the d\'ecalage pseudomonad restricts to this
$2$-subcategory; we denote its restriction by
$\Smon_{\Hcat}^{\mathrm{desc}}$.

Let $q\colon\dot\Delta\to\Delta$ be the corestriction of $\suc$,
so that
\begin{equation}
  j\circ q=\suc.
  \label{eq:q-corestriction}
\end{equation}
Write $\epsilon^j\colon j^*\circ\PsRan_j\Rightarrow\id$ for the counit
of the biadjunction, and write
\[
  \tau\colon(\epsilon^j\ast\id_{j^*})\circ(\id_{j^*}\ast\eta)
       \Longrightarrow\id_{j^*}
\]
for its invertible triangle modification.

\begin{proposition}[The descent--d\'ecalage comparison]
\label{prop:descent-successor-comparison}
Every shifted diagram $A\circ\suc$ has descent object $A(\mathsf1)$.
Consequently, $\Smon_{\Hcat}$ restricts to a pseudomonad
$\Smon_{\Hcat}^{\mathrm{desc}}$ on
$[\dot\Delta,\Hcat]_{PS}^{\mathrm{desc}}$.
For $A$ in this $2$-category, the composite
\begin{equation}
\begin{split}
  \chi_A\colon\Dmon_{\Hcat}(A)
  &\xrightarrow{\ \iota_{\Dmon_{\Hcat}(A)}\ }
  \Dmon_{\Hcat}(A)\circ\suc
  =(\Dmon_{\Hcat}(A)\circ j)\circ q\\
  &\xrightarrow{\ \epsilon^j_{A\circ j}\ast\id_q\ }
  (A\circ j)\circ q=\Smon_{\Hcat}^{\mathrm{desc}}(A)
\end{split}
\label{eq:chi-component}
\end{equation}
defines a pseudonatural transformation
$\chi\colon\Dmon_{\Hcat}\Rightarrow\Smon_{\Hcat}^{\mathrm{desc}}$.
Its compatibility with the points is the invertible modification
\begin{equation}
  \omega\colon\chi\circ\eta\Longrightarrow\iota,
  \label{eq:descent-successor-comparison}
\end{equation}
whose component is the following pasting:
\begin{equation}
\begin{aligned}
  \omega_A\colon\chi_A\circ\eta_A
  &=(\epsilon^j_{A\circ j}\ast\id_q)
       \circ\iota_{\Dmon_{\Hcat}(A)}\circ\eta_A\\
  &\xRightarrow{\ \cong\ }
    (\epsilon^j_{A\circ j}\ast\id_q)
       \circ(\eta_A\ast\id_{\suc})\circ\iota_A\\
  &\xRightarrow{\ (\tau_A\ast\id_q)\ast\id_{\iota_A}\ }
    \iota_A.
\end{aligned}
\label{eq:chi-point}
\end{equation}
Here the first cell is pseudonaturality of $\eta_A$ at $\delta$,
and the second is the triangle modification of the biadjunction.
Thus $(\chi,\omega)\colon(\Dmon_{\Hcat},\eta)
\to(\Smon_{\Hcat}^{\mathrm{desc}},\iota)$ is a pointed comparison.
\end{proposition}

\begin{proof}
In $\Cat$, all descent objects exist, so the displayed formulas give
the pointed comparison on the entire diagram $2$-category.
The multiplication of $\Smon_{\Cat}$ equips every shifted diagram
$B\circ\suc$ with a successor contraction;
\cref{thm:contraction-transfer,prop:effective-descent-pseudoalgebra}
therefore make it effective for descent. Applying this to
$B=\Hcat(X,-)\circ A$, pseudonaturally in $X$, shows that
$A(\mathsf1)$ represents the descent data of $A\circ\suc\circ j$.
Thus $\Smon_{\Hcat}$ restricts to
$[\dot\Delta,\Hcat]_{PS}^{\mathrm{desc}}$, where the same formulas
define $\chi$ and $\omega$.
\end{proof}

\begin{corollary}[Transfer to descent pseudoalgebras]
\label{cor:successor-descent-algebra}
The pointed comparison induces homomorphisms
\begin{equation}
  \Contr(\Smon_{\Hcat}^{\mathrm{desc}})
  \xrightarrow{\ \chi^*\ }
  \Contr(\Dmon_{\Hcat})
  \xrightarrow{\ E_{\Dmon_{\Hcat}}\ }
  \PsAlg(\Dmon_{\Hcat}),
  \label{eq:relative-descent-transfer}
\end{equation}
The transferred action and unit cell on a successor contraction
$(A,r,\varepsilon)$ are
\begin{equation}
  a_A=r\circ\chi_A,
  \qquad
  \bar a_A=\varepsilon\circ(\id_r\ast\omega_A).
  \label{eq:transferred-descent-action}
\end{equation}
In particular, $A$ is of effective $j$-descent.
\end{corollary}

\subsection{Postcomposition and absolute descent}

The algebraic structures just used are preserved under
postcomposition. This holds for the d\'ecalage pseudomonad without
any assumption about preservation of descent objects.

\begin{proposition}[Postcomposition of contractions and pseudoalgebras]
\label{prop:postcomposition-successor}
Every pseudofunctor $F\colon\Hcat\to\Kcat$ induces homomorphisms
\begin{equation}
\begin{aligned}
  F_*\colon\Contr(\Smon_{\Hcat})
       &\longrightarrow\Contr(\Smon_{\Kcat}),\\
  F_*\colon\PsAlg(\Smon_{\Hcat})
       &\longrightarrow\PsAlg(\Smon_{\Kcat}),
\end{aligned}
\label{eq:postcomposition-successor-algebras}
\end{equation}
compatible with the forgetful homomorphisms from pseudoalgebras to
contractions.
\end{proposition}

\begin{corollary}[Successor contractions give absolute descent]
\label{thm:absolute-descent}
Let $\Hcat$ be any $2$-category. If
$A\colon\dot\Delta\to\Hcat$ admits a successor contraction, then
$F\circ A$ is of effective descent for every pseudofunctor
$F\colon\Hcat\to\Kcat$. Hence $A$ is of absolute descent.
\end{corollary}

\begin{proof}
Postcomposition preserves successor contractions by
\cref{prop:postcomposition-successor}. Apply
\cref{cor:successor-descent-algebra} in $\Cat$ to
$\Kcat(X,-)\circ F\circ A$ for every $X\in\Kcat$; the resulting
equivalences exhibit $F(A(\mathsf0))$ as the descent object of
$F\circ A\circ j$.
\end{proof}

\subsection{Free d\'ecalage pseudoalgebras and absolute shifts}

The multiplication of the d\'ecalage pseudomonad supplies a canonical
source of the structures to which the transfer applies.

\begin{theorem}[Absolute left d\'ecalage]
\label{thm:absolute-shifts}
For every augmented pseudocosimplicial object
$A\colon\dot\Delta\to\Hcat$, the shift
\[
  A^+:=A\circ\suc=\Smon_{\Hcat}(A)
\]
is the carrier of the free $\Smon_{\Hcat}$-pseudoalgebra on $A$.
Its action is
\begin{equation}
  \id_A\ast\mu\colon A\circ\suc\circ\suc\longrightarrow A\circ\suc.
  \label{eq:free-successor-action}
\end{equation}
Its action and unit cell form a successor contraction.
Consequently, $A^+$ is of absolute descent.
\end{theorem}

\begin{corollary}[Absolute d\'ecalage]
\label{cor:absolute-decalage}
Let $G\colon\Delta\to\Hcat$ be pseudocosimplicial, and let
$q\colon\dot\Delta\to\Delta$ be the corestriction of
$\suc=\mathsf1+(-)$, so that $j\circ q=\suc$. Then
$G\circ q\colon\dot\Delta\to\Hcat$ carries a canonical
$\Smon_{\Hcat}$-pseudoalgebra structure and is of absolute descent.
Its action has components
\[
  G(\sigma_{n+1}^0)\colon G(\mathsf{n+2})\longrightarrow G(\mathsf{n+1}).
\]
\end{corollary}

\begin{proof}
The multiplication $\mu$ corestricts to a natural transformation
$\bar\mu\colon q\circ\suc\Rightarrow q$ in $\Delta$. The identities
\[
  \bar\mu\circ(\id_q\ast\delta)=\id_q,
  \qquad
  \bar\mu\circ(\bar\mu\ast\id_{\suc})
   =\bar\mu\circ(\id_q\ast\mu)
\]
are the unit and associativity laws for the action of $\suc$ on
$q$. Applying $G$ gives the pseudoalgebra action $\id_G\ast\bar\mu$ and
its coherence modifications. Forgetting multiplication and applying
\cref{thm:absolute-descent} proves the assertion.
\end{proof}

\begin{corollary}[Restriction of the d\'ecalage pseudomonad]
\label{cor:restricted-successor}
For arbitrary $\Hcat$, the endopseudofunctor $\Smon_{\Hcat}$ takes
every object of $[\dot\Delta,\Hcat]_{PS}$ into
$[\dot\Delta,\Hcat]_{PS}^{\mathrm{desc}}$. It therefore restricts to a
pseudomonad $\Smon_{\Hcat}^{\mathrm{desc}}$ on
$[\dot\Delta,\Hcat]_{PS}^{\mathrm{desc}}$.
The comparison of \cref{prop:descent-successor-comparison} is a
pointed comparison on this $2$-category
\[
  (\chi,\omega)\colon\Dmon_{\Hcat}
       \longrightarrow\Smon_{\Hcat}^{\mathrm{desc}}.
\]
The inclusion identifies the contractions and the pseudoalgebras of $\Smon_{\Hcat}^{\mathrm{desc}}$ with those of
$\Smon_{\Hcat}$.
\end{corollary}

Thus the full transfer may be displayed as
\begin{equation}
\begin{tikzcd}[column sep=7.0em,row sep=3.2em]
  \PsAlg(\Smon_{\Hcat}) \ar[r,"U_{\Smon_{\Hcat}}"]
  & \Contr(\Smon_{\Hcat}) \ar[d,"\chi^*"]
  \\
  \PsAlg(\Dmon_{\Hcat})
  & \Contr(\Dmon_{\Hcat})
     \ar[l,"E_{\Dmon_{\Hcat}}"',"\simeq"]
\end{tikzcd}
\label{eq:descent-algebra-transfer}
\end{equation}
All arrows are homomorphisms of $2$-categories. The lower arrow
supplies the uniquely determined multiplication coherence; it also
extends every pseudomorphism and $2$-cell of contractions.
Composing this diagram with the free-pseudoalgebra homomorphism
\[
  \mathsf{Fr}_{\Smon_{\Hcat}}\colon
  [\dot\Delta,\Hcat]_{PS}\longrightarrow\PsAlg(\Smon_{\Hcat})
\]
gives descent pseudoalgebras with underlying carriers $A\mapsto A^+$.
The absolute-descent theorem gives the additional assertion that
these carriers retain their descent universal property under every
pseudofunctor.

\section{Descent and Eilenberg--Moore completion}
\label{sec:BR}

We first isolate the implication from effective descent to monadicity,
which requires only a left adjoint to the augmentation
\cite{Nunes2021}. For the converse and the comparison of
factorizations, we apply the results of \cite{Nunes2018}, using
absolute d\'ecalage to supply the effective shifted row.

\subsection{Effective descent implies monadicity}
\label{subsec:effective-monadic}

The implication from effective descent to monadicity belongs to
a more general monadicity principle. The creation theorem of
\cite[Corollaries~2.6--2.7]{Nunes2021} shows that forgetful morphisms
from descent objects create absolute Kan extensions and, in
$\Cat$, absolute limits and colimits. Creation of absolute
coequalizers and Beck's monadicity theorem imply that a right adjoint
which factors, up to natural isomorphism, as an equivalence followed
by such a forgetful functor is monadic
\cite[Theorem~2.8]{Nunes2021}.

\Cref{prop:effective-implies-monadic} is the representable
form of this principle. This implication requires
no Beck--Chevalley condition; Beck--Chevalley enters the converse
and the comparison of the descent and Eilenberg--Moore
factorizations.

For a monad $T$ on $Y\in\Hcat$, its Eilenberg--Moore object $Y^T$
represents $Z\mapsto\operatorname{Alg}_{\Hcat(Z,T)}(\Hcat(Z,Y))$.
If $u\colon X\to Y$ has a left adjoint $u_!$, its counit induces
the comparison $K^T\colon X\to Y^T$ for $T=u\circ u_!$;
$u$ is \emph{monadic} when this comparison is an equivalence.

\begin{proposition}[Effective descent implies monadicity]
\label{prop:effective-implies-monadic}
Let $A\colon\dot\Delta\to\Hcat$ be of effective $j$-descent.
If $d=d_0^0\colon A(\mathsf 0)\to A(\mathsf 1)$ has a left adjoint $d_!$, then
$d$ is monadic. More precisely, $A(\mathsf 0)$, together with $d$ and the
action $\id_d\ast\epsilon$, is an Eilenberg--Moore object for $d\circ d_!$,
where $\epsilon\colon d_!\circ d\Rightarrow\id_{A(\mathsf 0)}$ is the counit.
\end{proposition}

\begin{proof}
Apply \cite[Theorem~2.8]{Nunes2021} to each $\Hcat(Z,d)$.
\end{proof}

In particular, only the converse implication from monadicity to
effective descent remains to be established under the additional
hypotheses of the formal B\'enabou--Roubaud theorem. Its comparison
of the two factorizations is stronger: it identifies their
intermediate objects even when their comparison morphisms are not
equivalences. 

\begin{example}[Monadicity alone does not imply effective descent]
\label{ex:monadic-not-effective}
This is the constant-diagram construction of
\cite[Remark~4.5 and Lemma~4.6]{Nunes2021}.
Let $U\colon\mathcal B\to\mathcal C$ be a monadic functor that is
not an equivalence, such as the forgetful functor from pointed sets
to sets. Augment the constant cosimplicial object at $\mathcal C$
by $U$: put $A(\mathsf 0)=\mathcal B$, $A(\mathsf n)=\mathcal C$ for $n\geq1$,
send every map between nonempty ordinals to the identity, and send
every map $\mathsf0\to\mathsf n$ with $n\geq1$ to $U$.
Normalization forces the descent isomorphism to be an identity,
so $\Desc(A\circ j)\simeq\mathcal C$ and $k_A$ identifies with $U$.
Thus $d_0^0$ is monadic, but $A$ is not of effective descent.
\end{example}

\subsection{The Beck--Chevalley comparison}

Let $\alpha\colon A\to B$ be a pseudonatural transformation of
augmented pseudocosimplicial objects, and suppose that
$\alpha_{n!}\dashv\alpha_n$ for every $n\geq0$. For a morphism
$v\colon\mathsf m\to\mathsf n$, orient its pseudonaturality cell as
$B(v)\circ\alpha_m\Rightarrow\alpha_n\circ A(v)$. Its left mate is
\begin{equation}
  \overline\alpha_v\colon
  \alpha_{n!}\circ B(v)\Longrightarrow A(v)\circ\alpha_{m!}.
  \label{eq:left-mate}
\end{equation}

\begin{definition}[Beck--Chevalley]
\label{def:BC}
The transformation $\alpha$ satisfies the \emph{Beck--Chevalley
condition} if all mates \eqref{eq:left-mate} are invertible.
\end{definition}

Under this condition, the adjunctions assemble into
$\alpha_!\dashv\alpha$ in $[\dot\Delta,\Hcat]_{PS}$
\cite[Lemma~7.3]{Nunes2018}.
Applied representably, \cite[Theorem~7.4]{Nunes2018} then gives,
without an ambient completeness assumption:
if $B$ is of effective descent and $\alpha_n$ is monadic for every
$n\geq1$, then
\begin{equation}
  A\text{ is of effective descent}
  \quad\Longleftrightarrow\quad
  \alpha_0\text{ is monadic}.
  \label{eq:general-BR}
\end{equation}
The hypotheses concern the diagram
\[
\begin{tikzcd}[column sep=3.8em,row sep=3.2em]
 A(\mathsf0) \ar[r,"d_0^0"] \ar[d,"\alpha_0"']
 & A(\mathsf1) \ar[r,shift left=1.4ex] \ar[r,shift right=1.4ex]
     \ar[d,"\alpha_1"']
 & A(\mathsf2) \ar[l]
     \ar[r,shift left=1.4ex] \ar[r] \ar[r,shift right=1.4ex]
     \ar[d,"\alpha_2"']
 & \cdots \\
 B(\mathsf0) \ar[r,"d_0^0"']
 & B(\mathsf1) \ar[r,shift left=1.4ex] \ar[r,shift right=1.4ex]
 & B(\mathsf2) \ar[l]
     \ar[r,shift left=1.4ex] \ar[r] \ar[r,shift right=1.4ex]
 & \cdots .
\end{tikzcd}
\]
For the d\'ecalage point $\id_A\ast\delta\colon A\to A\circ\suc$,
the corresponding mates are
\begin{equation}
  (d_n^0)_!\circ A(\suc(v))
  \Longrightarrow A(v)\circ(d_m^0)_!.
  \label{eq:successor-BC}
\end{equation}

\subsection{The formal B\'enabou--Roubaud theorem}

Let $A\colon\dot\Delta\to\Hcat$ be augmented pseudocosimplicial, and
put
\[
  d=d_0^0\colon A(\mathsf 0)\longrightarrow A(\mathsf 1).
\]
If $\Desc(A\circ j)$ exists, the canonical cone gives the descent
factorization
\begin{equation}
\begin{tikzcd}[column sep=6.0em,row sep=3.2em]
  A(\mathsf 0) \ar[r,"k_A"] \ar[dr,"d"']
  & \Desc(A\circ j) \ar[d,"u_A"]
  \\
  & A(\mathsf 1)
\end{tikzcd}
\label{eq:descent-factorization}
\end{equation}
Suppose $d_!\dashv d$, put
\begin{equation}
  T_A=d\circ d_!,
  \label{eq:TA}
\end{equation}
and assume that the Eilenberg--Moore object $A(\mathsf 1)^{T_A}$ exists. The
adjunction gives the factorization
\begin{equation}
\begin{tikzcd}[column sep=6.0em,row sep=3.2em]
  A(\mathsf 0) \ar[r,"K^{T_A}"] \ar[dr,"d"']
  & A(\mathsf 1)^{T_A} \ar[d,"U^{T_A}"]
  \\
  & A(\mathsf 1)
\end{tikzcd}
\label{eq:em-factorization}
\end{equation}

The following theorem identifies the two factorizations under
Beck--Chevalley. Its implication from effective descent to
monadicity is already covered by
\cref{prop:effective-implies-monadic}; the additional hypotheses
give the converse and the equivalence of the intermediate objects.

\begin{theorem}[Formal B\'enabou--Roubaud theorem]
\label{thm:formal-BR}
Let $\Hcat$ be any $2$-category and let
$A\colon\dot\Delta\to\Hcat$ be an augmented pseudocosimplicial
object. Suppose that
\[
  \id_A\ast\delta\colon A\longrightarrow A\circ\suc
\]
has componentwise left adjoints and satisfies the
Beck--Chevalley condition \eqref{eq:successor-BC}. Put
\[
  d=d_0^0,
  \qquad
  d_!=(d_0^0)_!,
  \qquad
  T_A=d\circ d_!.
\]
Whenever $\Desc(A\circ j)$ and $A(\mathsf 1)^{T_A}$ exist, there is an equivalence
\begin{equation}
  \Phi_A\colon\Desc(A\circ j)\simeq A(\mathsf 1)^{T_A}
  \label{eq:formal-BR-equivalence}
\end{equation}
and invertible $2$-cells
\begin{equation}
  U^{T_A}\circ\Phi_A\cong u_A,
  \qquad
  \Phi_A\circ k_A\cong K^{T_A}.
  \label{eq:formal-BR-compatibility}
\end{equation}
Thus the ladder
\begin{equation}
\begin{tikzcd}[column sep=4.2em,row sep=3.2em]
  A(\mathsf 0) \ar[r,"k_A"] \ar[d,equal]
  & \Desc(A\circ j) \ar[r,"u_A"] \ar[d,"\Phi_A"',"\simeq"]
  & A(\mathsf 1) \ar[d,equal]
  \\
  A(\mathsf 0) \ar[r,"K^{T_A}"']
  & A(\mathsf 1)^{T_A} \ar[r,"U^{T_A}"']
  & A(\mathsf 1)
\end{tikzcd}
\label{eq:BR-ladder}
\end{equation}
commutes up to the specified cells. In particular,
\[
  A\text{ is of effective }j\text{-descent}
  \quad\Longleftrightarrow\quad
  d_0^0\text{ is monadic}.
\]
\end{theorem}

\begin{proof}
Apply the comparison above to $\alpha=\id_A\ast\delta$.
The shifted object $A\circ\suc$ is of effective descent by
\cref{thm:absolute-shifts}. For $n\geq1$, the identity
$s_n^0\circ d_n^0\cong\id_{A(\mathsf n)}$ makes $d_n^0$ a
pseudosection, hence monadic by \cite[Lemma~8.4]{Nunes2018}.
The equivalence of factorizations is therefore
\cite[Theorem~8.5]{Nunes2018}, whose shift-effectiveness hypothesis
is supplied by \cref{thm:absolute-shifts}.
\end{proof}

Applied to $\Hcat(Z,-)\circ A$, this comparison identifies the
descent-data and algebra pseudofunctors, pseudonaturally in $Z$.
Thus either is representable precisely when the other is, giving
the following consequence.

\begin{corollary}[Local existence]
\label{cor:BR-local-existence}
Under the adjointness and Beck--Chevalley hypotheses of
\cref{thm:formal-BR}, either of $\Desc(A\circ j)$ and $A(\mathsf 1)^{T_A}$
exists if and only if the other does. In particular, existence of
both need not be imposed independently.
\end{corollary}

\section{Indexed \v{C}ech descent}
\label{sec:split}

We apply the preceding results to the indexed \v{C}ech objects of
Grothendieck descent. The formal B\'enabou--Roubaud theorem compares
their descent and Eilenberg--Moore factorizations. A section of the
base morphism gives a d\'ecalage pseudoalgebra. Its underlying
successor contraction then yields absolute descent, without any
adjointness assumption, and the Janelidze--Tholen characterization.

\subsection{The augmented \v{C}ech nerve}

Let $\C$ be a category with pullbacks and let $p\colon E\to B$. Put
\begin{equation}
  E^{[0]}=B,
  \qquad
  E^{[n]}=
  \underbrace{E\times_B\cdots\times_BE}_{n\text{ factors}}
  \quad(n\geq1).
  \label{eq:cech-powers}
\end{equation}
The augmented \v{C}ech nerve is the functor
\begin{equation}
  D_p\colon\dot\Delta^{\mathrm{op}}\longrightarrow\C,
  \qquad
  D_p(\mathsf n)=E^{[n]}.
  \label{eq:cech-nerve}
\end{equation}
For $v\colon\mathsf m\to\mathsf n$, the morphism
$D_p(v)\colon E^{[n]}\to E^{[m]}$ is characterized, for $m>0$, by
\[
  \pi_i^{(m)}\circ D_p(v)=\pi_{v(i)}^{(n)}
  \qquad(0\leq i<m).
\]
For $m=0<n$, it is the common morphism $E^{[n]}\to B$.

We also regard $D_p$ as a functor
$\dot\Delta\to\C^{\mathrm{op}}$. Given a pseudofunctor
\begin{equation}
  \Lx\colon\C^{\mathrm{op}}\longrightarrow\Hcat,
  \label{eq:indexed-pseudofunctor}
\end{equation}
put
\begin{equation}
  \Lx_p:=\Lx\circ D_p\colon\dot\Delta\longrightarrow\Hcat.
  \label{eq:indexed-cech-object}
\end{equation}
When it exists, write
\[
  \Desc_{\Lx}(p):=\Desc(\Lx_p\circ j).
\]
The comparison and forgetful $1$-cells are
\begin{equation}
  K_{\Lx}^{p}\colon\Lx(B)\longrightarrow\Desc_{\Lx}(p),
  \qquad
  U_{\Lx}^{p}\colon\Desc_{\Lx}(p)\longrightarrow\Lx(E).
  \label{eq:indexed-comparison-forgetful}
\end{equation}
They give the factorization
\begin{equation}
\begin{tikzcd}[column sep=6.0em,row sep=3.2em]
  \Lx(B) \ar[r,"K_{\Lx}^{p}"] \ar[dr,"\Lx(p)"']
  & \Desc_{\Lx}(p) \ar[d,"U_{\Lx}^{p}"]
  \\
  & \Lx(E)
\end{tikzcd}
\label{eq:indexed-descent-factorization}
\end{equation}
We say that $p$ is of \emph{effective $\Lx$-descent} if
$\Desc_{\Lx}(p)$ exists and $K_{\Lx}^{p}$ is an equivalence. For
$\Hcat=\Cat$, this is the classical category of descent data and its
usual factorization \cite{JanelidzeTholen1997,Nunes2018}. Taking
$\Lx=\C/(-)$ recovers the ordinary notion of an effective descent
morphism.

Before imposing Beck--Chevalley, we record the implication that
requires only a left adjoint to the single change-of-base morphism.

\begin{corollary}[Effective indexed descent implies monadicity]
\label{cor:indexed-effective-monadic}
If $p$ is of effective $\Lx$-descent and $\Lx(p)$ has a left
adjoint, then $\Lx(p)$ is monadic. No Beck--Chevalley condition is
required.
\end{corollary}

\begin{proof}
Apply \cref{prop:effective-implies-monadic} to $\Lx_p$, whose
augmentation is $\Lx(p)$.
\end{proof}

For $\Hcat=\Cat$, this is the implication established in
\cite[Theorem~4.7]{Nunes2021}. The preceding proof gives its
representable form for an arbitrary $2$-category.

\subsection{The indexed B\'enabou--Roubaud theorem}

Assume that every change-of-base $1$-cell $\Lx(f)$ has a left
adjoint. For a pullback square
\begin{equation}
\begin{tikzcd}[column sep=5.5em,row sep=3.2em]
  X' \ar[r,"v'"] \ar[d,"u'"']
     \ar[dr,phantom,"\lrcorner"{pos=0.14}]
  & X \ar[d,"u"]
  \\
  Y' \ar[r,"v"'] & Y
\end{tikzcd}
\label{eq:pullback-square}
\end{equation}
pseudofunctoriality gives an invertible $2$-cell
\[
  \Lx(v')\circ\Lx(u)\Longrightarrow\Lx(u')\circ\Lx(v).
\]
With the adjunctions on the two vertical change-of-base $1$-cells,
its left mate is
\begin{equation}
  \Lx(u')_!\circ\Lx(v')\Longrightarrow\Lx(v)\circ\Lx(u)_!.
  \label{eq:indexed-BC-mate}
\end{equation}
We say that $\Lx$ satisfies \emph{Beck--Chevalley} if this mate is
invertible for every pullback square. Equivalently one may mate with
respect to the horizontal sides: transposing the same pullback square
gives that convention.

Fix $p$ and put $A=\Lx_p$. For $n\geq0$, let
\[
  f_n=D_p(\delta_n)\colon E^{[n+1]}\longrightarrow E^{[n]}.
\]
For $v\colon\mathsf m\to\mathsf n$, naturality of $\delta$ gives
the square
\begin{equation}
\begin{tikzcd}[column sep=6.0em,row sep=3.4em]
  E^{[n+1]} \ar[r,"D_p(\suc(v))"] \ar[d,"f_n"']
    \ar[dr,phantom,"\lrcorner"{pos=0.14}]
  & E^{[m+1]} \ar[d,"f_m"]
  \\
  E^{[n]} \ar[r,"D_p(v)"'] & E^{[m]}
\end{tikzcd}
\label{eq:cech-BC-square}
\end{equation}
This is a pullback. Indeed, under
$E^{[n+1]}\cong E\times_BE^{[n]}$, it is obtained by pulling back
$E\times_BE^{[m]}\to E^{[m]}$ along $D_p(v)$. After applying $\Lx$,
its pseudofunctoriality cell is precisely the pseudonaturality cell of
$\id_A\ast\delta$. Hence Beck--Chevalley for $\Lx$ implies the condition
\eqref{eq:successor-BC} for $A$.

\begin{theorem}[Two-dimensional B\'enabou--Roubaud]
\label{thm:two-dimensional-BR}
Let $\C$ be a category with pullbacks, let $\Hcat$ be any
$2$-category, and let
$\Lx\colon\C^{\mathrm{op}}\to\Hcat$ be a pseudofunctor. Suppose that
every $\Lx(f)$ has a left adjoint and that $\Lx$ satisfies
Beck--Chevalley. For $p\colon E\to B$, put
\[
  T_p=\Lx(p)\circ\Lx(p)_!.
\]
The objects $\Desc_{\Lx}(p)$ and $\Lx(E)^{T_p}$ exist simultaneously.
Whenever they exist, there is an equivalence
\begin{equation}
  \Phi_p\colon
  \Desc_{\Lx}(p)\simeq\Lx(E)^{T_p}
  \label{eq:indexed-BR-equivalence}
\end{equation}
with invertible $2$-cells
\[
  U^{T_p}\circ\Phi_p\cong U_{\Lx}^{p},
  \qquad
  \Phi_p\circ K_{\Lx}^{p}\cong K^{T_p}.
\]
Consequently,
\[
  p\text{ is of effective }\Lx\text{-descent}
  \quad\Longleftrightarrow\quad
  \Lx(p)\text{ is monadic}.
\]
\end{theorem}

\begin{proof}
For $A=\Lx_p$, the component $d_0^0$ is $\Lx(p)$ and the induced
monad is $T_p$. The pullback calculation
\eqref{eq:cech-BC-square} supplies Beck--Chevalley for $\id_A\ast\delta$.
Apply \cref{thm:formal-BR,cor:BR-local-existence}.
\end{proof}

\begin{corollary}[B\'enabou--Roubaud]
\label{cor:classical-BR}
Let $\C$ be a category with pullbacks and let
$\Lx\colon\C^{\mathrm{op}}\to\Cat$ be an indexed category whose
change-of-base functors have left adjoints and satisfy
Beck--Chevalley. Then, for every $p\colon E\to B$, the descent
factorization of $\Lx(p)$ is equivalent to its Eilenberg--Moore
factorization. In particular,
\[
  p\text{ is of effective }\Lx\text{-descent}
  \quad\Longleftrightarrow\quad
  \Lx(p)\text{ is monadic}.
\]
\end{corollary}

\begin{proof}
The required descent categories and Eilenberg--Moore categories exist
in $\Cat$, so this is \cref{thm:two-dimensional-BR} with
$\Hcat=\Cat$.
\end{proof}

\subsection{The \v{C}ech contraction of a split epimorphism}

Suppose now that $p\colon E\to B$ has a section
\begin{equation}
  s\colon B\longrightarrow E,
  \qquad
  p\circ s=\id_B.
  \label{eq:section-of-p}
\end{equation}
Define
\begin{equation}
  h_n\colon E^{[n]}\longrightarrow E^{[n+1]}
  \label{eq:hn}
\end{equation}
by $h_0=s$ and, for $n\geq1$, by
\begin{equation}
  \pi_0^{(n+1)}\circ h_n=s\circ p\circ\pi_0^{(n)},
  \qquad
  \pi_{i+1}^{(n+1)}\circ h_n=\pi_i^{(n)}
  \quad(0\leq i<n).
  \label{eq:hn-projections}
\end{equation}
In element notation, this inserts the section in the first
coordinate:
\[
  h_n(e_0,\ldots,e_{n-1})
  =\bigl((s\circ p)(e_0),e_0,\ldots,e_{n-1}\bigr).
\]

\begin{lemma}[The \v{C}ech contraction]
\label{lem:cech-contraction}
For every $v\colon\mathsf m\to\mathsf n$ in $\dot\Delta$,
\begin{align}
  D_p(\suc(v))\circ h_n&=h_m\circ D_p(v),
  \label{eq:cech-contraction-naturality}\\
  D_p(\delta_n)\circ h_n&=\id_{E^{[n]}},
  \label{eq:cech-contraction-retraction}\\
  h_{n+1}\circ h_n&=D_p(\mu_n)\circ h_n.
  \label{eq:cech-contraction-multiplication}
\end{align}
Consequently, regarding $D_p$ as an augmented cosimplicial object in
$\C^{\mathrm{op}}$, the opposites of the $h_n$ form a strict
$\Smon_{\C^{\mathrm{op}}}$-algebra, whose underlying contraction
is the successor contraction.
\end{lemma}

\begin{proof}
The first two identities follow by comparing projections, using
$p\circ s=\id_B$. For the third, both sides insert the section twice.
\end{proof}

In low degrees,
\[
  h_0=s,
  \qquad
  h_1=\langle s\circ p,\id_E\rangle
  \colon E\to E\times_BE.
\]
The opposites of these equations give the ordinary split-equalizer
identities. The higher maps $h_n$ give the full successor algebra.

Postcomposition with $\Lx$ transports this algebra to a
$\Smon_{\Hcat}$-pseudoalgebra on $A=\Lx_p$. Its action is
\begin{equation}
  r\colon A\circ\suc\longrightarrow A,
  \qquad r_n=\Lx(h_n).
  \label{eq:indexed-successor-action}
\end{equation}
The naturality squares \eqref{eq:cech-contraction-naturality},
together with the compositors of $\Lx$, supply its pseudonaturality
cells. Its unit and multiplication cells have components
\begin{equation}
\begin{aligned}
  r_n\circ d_n^0
  &=\Lx(h_n)\circ\Lx(D_p(\delta_n))
    \cong\Lx(D_p(\delta_n)\circ h_n)
    \cong\id_{A(\mathsf n)},\\
  r_n\circ r_{n+1}
  &=\Lx(h_n)\circ\Lx(h_{n+1})
    \cong\Lx(h_{n+1}\circ h_n)\\
  &=\Lx(D_p(\mu_n)\circ h_n)
    \cong\Lx(h_n)\circ\Lx(D_p(\mu_n)).
\end{aligned}
\label{eq:indexed-successor-coherence}
\end{equation}
The last expression is the degree-$n$ component of $r\circ m_A$,
where $m_A=\id_A\ast\mu$ is the d\'ecalage multiplication. Thus the second
line gives the modification $r\circ\Smon(r)\Rightarrow r\circ m_A$.
The strict algebra laws just proved and the pseudofunctor coherence
axioms for $\Lx$ give the pseudoalgebra unit and associativity
axioms. In particular, the multiplication coherence is part of the
induced structure, rather than an additional choice.

\begin{theorem}[Absolute indexed descent for split epimorphisms]
\label{thm:split-absolute}
Let $\C$ be a category with pullbacks and let $p\colon E\to B$ be a
split epimorphism. For every $2$-category $\Hcat$ and every
pseudofunctor $\Lx\colon\C^{\mathrm{op}}\to\Hcat$, a
section of $p$ equips $\Lx_p$ with the $\Smon_{\Hcat}$-pseudoalgebra
\eqref{eq:indexed-successor-action}--\eqref{eq:indexed-successor-coherence}.
Its underlying successor contraction induces a descent pseudoalgebra,
and its augmented cone is an absolute descent cone with vertex
$\Lx(B)$. In particular, $p$ is of effective $\Lx$-descent.
\end{theorem}

\begin{proof}
The preceding construction gives the d\'ecalage pseudoalgebra.
Forgetting its multiplication cell and applying
\cref{thm:absolute-descent} proves existence and absoluteness of the
descent cone. Thus $\Lx_p$ belongs to the domain of descent completion
in \cref{prop:restricted-descent-completion}. Its descent action is
$r\circ\chi_{\Lx_p}\colon\Dmon(\Lx_p)\to\Lx_p$, with unit induced by the
pointed comparison and multiplication supplied by
\cref{prop:idempotent-upgrade}. This is the structure obtained by
the transfer of contractions.
\end{proof}

\subsection{The Janelidze--Tholen characterization}

A morphism $p$ is \emph{of absolute effective descent} if it is of
effective $\Lx$-descent for every indexed category
$\Lx\colon\C^{\mathrm{op}}\to\Cat$
\cite{JanelidzeTholen1997}. Effective $\Lx$-descent makes $\Lx(p)$
conservative, since the forgetful functor from descent data reflects
isomorphisms.

For a set $S$, put $\two^S=[S_{\mathrm{disc}},\two]$, where
$\two=(0\to1)$.
\begin{lemma}
\label{lem:restriction-surjective}
For a function $v\colon S\to T$, restriction
$v^*\colon\two^T\to\two^S$ is conservative if and only if $v$ is
surjective.
\end{lemma}
\begin{proof}
Surjectivity detects invertibility at every coordinate. If a coordinate
is missing, the transformation $0\to1$ there and identities elsewhere
is noninvertible but restricts to an identity.
\end{proof}

Fix $B\in\C$ and define the indexed category
\begin{equation}
  \Lx_B(X)=\two^{\C(B,X)},\qquad
  \Lx_B(f)=\C(B,f)^*.
  \label{eq:detecting-indexed-category}
\end{equation}
\begin{proposition}[One indexed category detects splitting]
\label{prop:detects-splitting}
If $p\colon E\to B$ is of effective $\Lx_B$-descent, then $p$ is a
split epimorphism.
\end{proposition}
\begin{proof}
Conservativity of $\Lx_B(p)$ makes $\C(B,p)$ surjective by
\cref{lem:restriction-surjective}. A preimage of $\id_B$ is a section
of $p$.
\end{proof}

\begin{theorem}[Janelidze--Tholen]
\label{thm:JT}
Let $\C$ be a category with pullbacks. For a morphism
$p\colon E\to B$, the following conditions are equivalent:
\begin{enumerate}[label=\textup{(\roman*)},leftmargin=2.7em]
  \item $p$ is of absolute effective descent;
  \item $p$ is a split epimorphism.
\end{enumerate}
\end{theorem}

\begin{proof}
If $p$ is of absolute effective descent, it is of effective
$\Lx_B$-descent and therefore splits by
\cref{prop:detects-splitting}. Conversely, if $p$ splits,
\cref{thm:split-absolute} with $\Hcat=\Cat$ shows that it is of
effective $\Lx$-descent for every indexed category $\Lx$.
\end{proof}

\section{Co-duality and co-Eilenberg--Moore objects}
\label{sec:codual}

Since the preceding results hold in an arbitrary $2$-category, they
also apply in $\Hcat^{\mathrm{co}}$. Reversing $2$-cells turns right
adjoints into left adjoints and comonads into monads. Conical bilimits,
and hence descent objects, are unchanged: inverting cone constraints
identifies the cone categories with their opposites.

For an adjunction $d\dashv d_*$, with unit $\eta$ and counit
$\epsilon$, write
\begin{equation}
  Q_A=d\circ d_*,\qquad
  \epsilon\colon Q_A\Rightarrow\id_{A(\mathsf1)},\qquad
  \id_d\ast\eta\ast\id_{d_*}\colon Q_A\Rightarrow Q_A\circ Q_A,
  \label{eq:codual-comonad}
\end{equation}
where $d=d_0^0\colon A(\mathsf0)\to A(\mathsf1)$.
The \emph{co-Eilenberg--Moore object} $A(\mathsf1)^{Q_A}$ represents
$Q_A$-coalgebras; equivalently, it is the Eilenberg--Moore object of
the corresponding monad in $\Hcat^{\mathrm{co}}$. The coaction
$\id_d\ast\eta$ induces the factorization
\begin{equation}
  A(\mathsf0)\xrightarrow{K^{Q_A}}A(\mathsf1)^{Q_A}
       \xrightarrow{U^{Q_A}}A(\mathsf1),\qquad
  U^{Q_A}\circ K^{Q_A}\cong d.
  \label{eq:coem-factorization}
\end{equation}
The $1$-cell $d$ is \emph{comonadic} when $K^{Q_A}$ is an equivalence.

\begin{corollary}[Effective descent implies comonadicity]
\label{prop:effective-implies-comonadic}
If $A$ is of effective descent and $d_0^0$ has a right adjoint,
then $d_0^0$ is comonadic. No Beck--Chevalley condition is required.
\end{corollary}
\begin{proof}
Apply \cref{prop:effective-implies-monadic} in $\Hcat^{\mathrm{co}}$.
\end{proof}

Suppose each $d_n^0$ has a right adjoint $r_n$.
The \emph{right Beck--Chevalley condition} requires the right mates
of the inverses of the pseudonaturality cells
$\gamma_v\colon A(\suc(v))\circ d_m^0\Rightarrow d_n^0\circ A(v)$
to be invertible:
\begin{equation}
  \gamma_v^*\colon A(v)\circ r_m\Longrightarrow
  r_n\circ A(\suc(v))\qquad(v\colon\mathsf m\to\mathsf n).
  \label{eq:right-successor-BC}
\end{equation}
This is the left Beck--Chevalley condition in $\Hcat^{\mathrm{co}}$.

\begin{corollary}[Co-dual B\'enabou--Roubaud]
\label{thm:codual-BR}
Let $A\colon\dot\Delta\to\Hcat$ be an augmented pseudocosimplicial
object. If each $d_n^0$ has a right adjoint and the right
Beck--Chevalley mates \eqref{eq:right-successor-BC} are invertible,
then $\Desc(A\circ j)$ and
$A(\mathsf1)^{Q_A}$ exist simultaneously. When they exist, there is an
equivalence
\begin{equation}
  \Phi_A^{\mathrm{co}}\colon\Desc(A\circ j)\simeq A(\mathsf1)^{Q_A},
  \label{eq:codual-BR-equivalence}
\end{equation}
with $U^{Q_A}\circ\Phi_A^{\mathrm{co}}\cong u_A$ and
$\Phi_A^{\mathrm{co}}\circ k_A\cong K^{Q_A}$.
In particular, $A$ is of effective descent if and only if $d_0^0$
is comonadic.
\end{corollary}
\begin{proof}
Apply \cref{thm:formal-BR,cor:BR-local-existence} in
$\Hcat^{\mathrm{co}}$.
\end{proof}

For an indexed pseudofunctor $\Lx\colon\C^{\mathrm{op}}\to\Hcat$
whose change-of-base $1$-cells have right adjoints, with $\C$ admitting
pullbacks, the corresponding
mates for \eqref{eq:pullback-square} are
\begin{equation}
  \Lx(v)\circ\Lx(u)_*\Longrightarrow
  \Lx(u')_*\circ\Lx(v').
  \label{eq:right-indexed-BC}
\end{equation}
They are the right mates of the inverses of the pseudofunctoriality
cells used in \eqref{eq:indexed-BC-mate}.

\begin{corollary}[Comonadic indexed descent]
\label{cor:comonadic-indexed-descent}
If the mates \eqref{eq:right-indexed-BC} are invertible for every
pullback square, then the descent and co-Eilenberg--Moore
factorizations of $\Lx(p)$ are equivalent whenever either exists.
Thus $p$ is of effective $\Lx$-descent if and only if $\Lx(p)$ is
comonadic. For $\Hcat=\Cat$, both objects always exist.
\end{corollary}
\begin{proof}
Apply \cref{thm:two-dimensional-BR} to the indexed pseudofunctor
with values in $\Hcat^{\mathrm{co}}$.
\end{proof}

\paragraph{Acknowledgements.}
I thank Walter Tholen for valuable discussions on descent theory during
my stay at the Fields Institute in August 2023.

\paragraph{Funding.}
This work was supported by the Centre for Mathematics of the
University of Coimbra (CMUC), funded by the Portuguese Foundation for
Science and Technology (FCT), under projects UID/00324/2025 and
UID/PRR/00324/2025.
I also acknowledge support from the Fields Institute for Research in
Mathematical Sciences through a Fields Research Fellowship in 2023.

\end{document}